\documentclass[10pt]{amsart}
\usepackage{amsmath,mathrsfs}
\usepackage{amssymb}
\usepackage{latexsym}
\usepackage{amsfonts}
\usepackage{graphicx}
\usepackage{psfrag}

\newtheorem{Pa}{Paper}[section]
\newtheorem{theorem}[Pa]{{\bf Theorem}}
\newtheorem{lemma}[Pa]{{\bf Lemma}}
\newtheorem{corollary}[Pa]{{\bf Corollary}}

\newtheorem{proposition}[Pa]{{\bf Proposition}}
\newtheorem{observation}[Pa]{{\bf Observation}}

\newtheorem{example}[Pa]{{\bf Example}}

\def\C{\mathbb C}
\def\R{{\mathbb R}}

\def\@tempb{saamsart}
\usepackage{hyperref}
\begin{document}
\bibliographystyle{plain}
\thispagestyle{empty}
\title[
realization of even generalized
positive and odd rational functions]
{state space realization of even
generalized positive and odd rational function.\\
Applications to static output feedback}
\author[D. Alpay]{Daniel Alpay}
\address{(DA) Department of mathematics,
Ben-Gurion University of the Negev, P.O. Box 653, Beer-Sheva
84105, Israel} \email{dany@math.bgu.ac.il}
\author[I. Lewkowicz]{Izchak ~Lewkowicz}
\address{(IL) Department of electrical
engineering, Ben-Gurion University of the Negev, P.O. Box 653,
Beer-Sheva 84105, Israel} \email{izchak@ee.bgu.ac.il}

\thanks{
This research was supported in part by the Bi-national Science
Foundation grant no. 2010117.
D. Alpay thanks the
Earl Katz family for endowing the chair
which supported his research.}

\def\squarebox#1{\hbox to #1{\hfill\vbox to #1{\vfill}}}

\maketitle
\begin{abstract}\quad
We here specialize the well known Positive Real Lemma (also
known as the Kalman-Yakubovich-Popov Lemma) to complex
matrix-valued rational functions, (i) generalized positive
even and (ii) odd. On the way we characterize the (non)
minimality of realization of arbitrary systems through
(i) the corresponding state matrix and (ii) moving the
poles by applying static output feedback. We then explore
the application of static output feedback to both
generalized positive even and to odd functions.
\end{abstract}

\section{Introduction}
\setcounter{equation}{0}

For a half of a century, the Positive Real Lemma (also known
as the Kalman-Yakubovich-Popov Lemma) has been recognized as
a fundamental result in System Theory. We here exploit it to
study two classes of rational functions, (i) generalized
positive even and (ii) odd.
\vskip 0.2cm

Let $\C_+$ and $\C_-$ be the open right and left halves of
the complex plane respectively, and $~{\mathbb P}_k,~
(\overline{\mathbb P}_k)~$ be the set of all $~k\times k~$
positive definite (semidefinite) matrices. Recall that a
$p\times p$-valued function $~F(s)$, analytic in
$\mathbb C_+$ is said to be {\em positive} if
\begin{equation}\label{eq:PosFunc}
\begin{matrix}
F(s)+F^*(s)
\in\overline{\mathbb P}_p &~&s\in\C_+~.\end{matrix}
\end{equation}
The study of rational positive functions, denoted by
${\mathcal P}$, has been motivated from the 1920's by
(lumped) electrical networks theory, see e.g.
\cite{AV}, \cite{Be}. From the 1960's positive
functions also appeared in books on absolute
stability theory, see e.g. \cite{NT}, \cite{Po}.

A $~p\times p$-valued function of bounded type in
$\mathbb C_+$ (i.e. a quotient of two functions
analytic and bounded in $\mathbb C_+$) is called
{\em generalized positive} $~\mathcal{GP}~$ if
\begin{equation}\label{DefGP}
F(i\omega)+F^*(i\omega)\in\overline{\mathbb P}_p
\quad a.e.\quad\quad\omega\in\R,
\end{equation}
where $F(i\omega)$ denotes the non-tangential
limit\begin{footnote}{This limit exists almost
everywhere on $i\R$ because $F$ is assumed of bounded
type in $\C_+$, see e.g. \cite{Du}.}\end{footnote}
of $F$ at the point $i\omega$.
\vskip 0.2cm

It is interesting to note that both the function set
$\mathcal{GP}$ and its subset $\mathcal P$, are closed
under positive scaling, sum and inversion (when the
given function has a non-identically vanishing
determinant). We emphasize that for simplicity we adhere
hereafter to the rational case. We shall find it
convenient to use the notation
\[
F^{\#}(s):=F^*(-s^*).
\]
Recall that a matrix valued function $F(s)$ is called
{\em even}~ if
\begin{equation}\label{even}
F(s)=F^{\#}(s).
\end{equation}
We shall denote by ${\mathcal Even}$ the set of even functions.
Having $F\in{\mathcal Even}$ in particular implies that its
zeroes and poles are symmetric with respect to the imaginary
axis. Moreover, on the imaginary axis this $F$
is Hermitian, i.e.
\[
F(i\omega)=F^*(\omega)\quad a.e. \quad\quad\omega\in\R.
\]
We shall say that a $p\times p$-valued
function\begin{footnote}{To ease reading, hereafter
$F(s)$ denotes an arbitrary rational function and
$\Psi(s)$ is from the subset of $\mathcal{GP}$.}
\end{footnote} $\Psi(s)$ is~ {\em even generalized
positive}, denoted by $\mathcal{GPE}$, if it satisfies
both \eqref{DefGP} and \eqref{even}, i.e.
\begin{equation}\label{DefGPE}
\Psi(i\omega)\in\overline{\mathbb P}_p
\quad a.e.\quad\quad\omega\in\R.
\end{equation}
Namely, $\mathcal{GPE}=\mathcal{GP}\bigcap{\mathcal Even}$.
Note having that $\Psi\in\mathcal{GP}$ is equivalent to
\mbox{$\left(\Psi+\Psi^{\#}\right)\in\mathcal{GPE}$}.
Scalar rational $\mathcal{GPE}$ functions were recently
studied and then applied to Nevanlinna-Pick interpolation
in \cite[Section 5]{AL3}.
\vskip 0.2cm

In analogy to ${\mathcal Even}$ functions in \eqref{even},
we shall say that a function $F(s)$ is ~{\em odd}, denoted
by $F\in{\mathcal Odd}$, if
\begin{equation}\label{odd}
F(s)=-F^{\#}(s).
\end{equation}
This implies that on the imaginary axis {\em Odd}
functions are skew-Hermitian,
\begin{equation}\label{eq:SkewHerm}
F(i\omega)=-F^*(i\omega)\quad a.e. \quad\quad\omega\in\R.
\end{equation}

Recall also that a matrix valued function $F(s)$ can
always be partitioned to its {\em even}~ and {\em odd}
parts, i.e.
\begin{center}
$F(s)=F_{\rm even}(s)+F_{\rm odd}(s)\quad\quad\quad
F_{\rm even}:=${\mbox{\tiny$\frac{1}{2}$}}$\left(
F+F^{\#}\right)\quad\quad\quad F_{\rm odd}:=$
{\mbox{\tiny$\frac{1}{2}$}}$\left(F-F^{\#}\right).$
\end{center}
Clearly, $F_{\rm even}\in{\mathcal Even}$ and
$F_{\rm odd}\in{\mathcal Odd}$. From \eqref{DefGP} and
\eqref{eq:SkewHerm} it follows that in fact
${\mathcal Odd}\subset\mathcal{GP}$. Scalar rational
${\mathcal Odd}$ functions were recently studied in
\cite[Section 4]{AL3}. To further motivate the study of
${\mathcal Odd}$ functions, recall that the classical
Nevanlinna-Pick framework it was shown in \cite{YS} that
if there exists an interpolating function within
$\mathcal{P}$, without loss of generality (but
compromising the minimal degree) there exists an
interpolating function within
\mbox{$\mathcal{PO}:=\mathcal{P}\bigcap{\mathcal Odd}$.}
For an application of this observation see
\cite[Corollary 5.2.2]{CL4}.
\vskip 0.2cm

Specializing the above even-odd partitioning to
$\Psi\in\mathcal{GP}$ one obtains
\[
\Psi(s)=\Psi_{\rm even}(s)+\Psi_{\rm odd}(s)
\quad\quad\quad\quad
\Psi_{\rm even}\in\mathcal{GPE}\quad\quad
\Psi_{\rm odd}\in{\mathcal Odd}.
\]
As already mentioned, scalar, rational $\mathcal{GPE}$ and
${\mathcal Odd}$ functions were recently studied in
\cite[Section 5]{AL3} and \cite[Section 4]{AL3}, respectively.
\vskip 0.2cm

We now recall in the concept of state space realization. Let
$F(s)$ be a $p\times p$-valued rational function analytic at
infinity, i.e.  $\lim\limits_{s~\rightarrow~\infty}~F(s)~$
exists. Denote by $q$ the McMillan degree of $F$. Namely,
$F$ admits a state space realization
\begin{equation}\label{StateSpace}
\begin{matrix}
F(s)=C(sI-A)^{-1}B+D&~&~&
L:=\left(\begin{smallmatrix}A&~B\\C&~D
\end{smallmatrix}\right)\end{matrix}
\end{equation}
with $~A\in\C^{n\times n}$, $n\geq q$, $~B, C^*\in\C^{n\times
p}~$ and $~D\in\C^{p\times p}$, namely, $~L\in\C^{(n+p)
\times(n+p)}$. If the McMillan degree of $F(s)~$ satisfies
$~q=n$, the realization is called minimal.
\vskip 0.2cm

The aim of this work is to characterize the realization of
rational $\mathcal{GPE}$ and ${\mathcal Odd}$ functions
and to give an application to static output feedback.
\vskip 0.2cm

Related works on realization for the non-rational case exist,
see e.g. \cite{DLLS2} and more recently \cite{ABDT}.
Restricting the discussion to rational functions, enables
to offer simple explicit formulas. On realization, our main
result is as follows.

\begin{theorem}\label{main}
Let $\Psi(s)$ be a $p\times p$-valued rational function so
that $~\lim\limits_{s~\rightarrow~\infty}~\Psi(s)~$ exists
and let $L$ be the associated realization matrix,
\[
\begin{matrix}
\Psi(s)=C(sI-A)^{-1}B+D&~&~&
L:=\left(\begin{smallmatrix}A&~B\\C&~D
\end{smallmatrix}\right)\end{matrix},
\]
see \eqref{StateSpace}. Let $q$ denotes the
McMillan degree.
\vskip 0.2cm

{\bf A}. Let $L\in\C^{(2n+p)\times(2n+p)}$ with $q$ even,
$2n\geq q$. The following are equivalent.
\begin{itemize}

\item[(i)~~~]{} $\Psi\in\mathcal{GPE}$.

\item[(ii)~~]{} There exist $~2n\times 2n$ Hermitian
non-singular matrices $H_1~$ and $~H_2~$ so that,
\[
HL+L^*H\in\overline{\mathbb{P}}_{2n+p}\quad\quad\quad
ML+(ML)^*=0_{2n+p},
\]
where
\[
H:={\rm diag}\{H_1,~I_p\}\quad\quad\quad
M:={\rm diag}\{H_2,~iI_p\}.
\]
\item[(iii)~]{} There exist realization matrices $L$
whose sub-blocks are
\begin{equation}\label{eq:L}
A=\left(\begin{smallmatrix}\hat{A}&~\hat{B}\hat{B}^*\\
0&-\hat{A}^*\end{smallmatrix}\right)\quad\quad
B=\left(\begin{smallmatrix}\hat{B}\hat{D}^*\\
-\hat{C}^*\end{smallmatrix}\right)\quad\quad
C=(\hat{C}\quad\hat{D}\hat{B}^*)\quad\quad
D=\hat{D}\hat{D}^*
\end{equation}
where $\hat{A}$ is $n\times n$, $\hat{B}, \hat{C}^*$ are
$n\times p$ and $\hat{D}$ is $p\times p$. Furtheremore without
loss of generality one can take spectrum $\hat{A}$ to be in
$\overline{\mathbb{C}_+}$.

\noindent
Moreover, if the realization is minimal, i.e. $2n=q$,
up to similarity, $L$ is of the above form.
\end{itemize}
\vskip 0.2cm

{\bf B}. Let $L\in\C^{(n+p)\times(n+p)}$ with $n\geq q$ and let
$\nu$ be so that $A$, the $n\times n$ part of $L$ has at
most $\nu$ eigenvalues in $\C_-$ and at most $n-\nu$
eigenvalues in $\C_+.$

Consider the following statements.
\begin{itemize}
\item[(i)~~~]{}$\Psi\in{\mathcal Odd}$

\item[(ii)~~]{}There exists a
non-singular $n\times n$ Hermitian $\hat{H}$ s.t.
\mbox{$H:={\rm diag}\{\hat{H}~,~I_p\}$} satisfies
\[
HL+L^*H=0_{n+p}~.
\]

\item[(iii)~]{} Up to similarity the sub-blocks of the
realization matrix $L$ are
\begin{equation}\label{RealOdd}
A=\left(\begin{smallmatrix} T_1&\tilde{A}\\
\tilde{A}^*&T_2\end{smallmatrix}\right)\quad\quad
B=\left(\begin{smallmatrix}B_1\\
B_2\end{smallmatrix}\right)\quad\quad
C=(B_1^*~~~-B_2^*)\quad\quad D=T_3
\end{equation}
with $T_1, T_2, T_3$ skew-Hermitian ($T_j=-T_j^*$) of
dimensions $\nu\times\nu$, \mbox{$(n-\nu)\times(n-\nu)$}
and $p\times p$, respectively.
\end{itemize}
Then $(ii)$ and $(iii)$ are equivalent and in turn imply
$(i)$.

If in addition the realization is minimal, i.e. $q=n$,
the the converse is true as well.
\end{theorem}
\vskip 0.2cm

The result of part {\bf A} can be compared with
\cite[Theorem 10.2]{BGKR} where a state space realization
of $\mathcal{GPE}$ functions is given. While the
development in \cite{BGKR} is self contained, we
here rely on two classical results: Lemma \ref{gpLemma}
and Theorem \ref{Pr:GPEchar} below. This allows us to
obtain the above Theorem \ref{main} {\bf A} which is more
general than \cite[Theorem 10.2]{BGKR} in the following
sense, (i) minimality of the realization is not assumed
and (ii) $D=\lim\limits_{s\rightarrow\infty}\Psi(s)$
can be of any rank (including zero). Moreover, the explicit
formulation in \eqref{eq:L} and \eqref{RealOdd} enables
us to directly address applications like static output
feedback. For example, recall that a $\mathcal{GPE}$
function which is analytic on $i\R$ admits a spectral
factorization $GG^{\#}$ with $G(s)$ analytic in
$\overline{\C_+}$, see e.g. \cite[Chapter 9]{BGKR},
\cite{Fu}, \cite[Section 19.3]{LR} and \cite{Ran}.
We here state a result, where some of the details will
be clarified in the sequel.

\begin{proposition}\label{FeedBackSpectFact}
Let $\Psi(s)$ be a $p\times p$-valued rational
$\mathcal{GPE}$ function, which is not analytic on $i\R$.
Assume that $\lim\limits_{s\rightarrow\infty}\Psi(s)=0$
and let the state space realization be as in \eqref{eq:L}
with $\hat{D}=0$.
\vskip 0.2cm

There exists a static output feedback gain $K$,
$-K\in\overline{\mathbb P}_p$ so that the closed loop
system \mbox{$(I-\Psi{K})^{-1}\Psi$} is analytic on $i\R$,
if and only if for all $r\in\R$ the two following matrices
\[
\left(\begin{smallmatrix}
\hat{A}-irI_n&\vdots&\hat{B}\end{smallmatrix}\right)
\quad\quad{and}\quad\quad
\left(\begin{smallmatrix}
\hat{A}-irI_n\\ \hat{C}\end{smallmatrix}\right)
\]
are of full rank.
\end{proposition}
\vskip 0.2cm

Clearly, it is suffices to check the conditions only for
all $ir\in{\rm spect}(\hat{A})$.
\vskip 0.2cm

The outline of the paper is as follows. Section \ref{Sec:2}
is devoted to providing a perspective on relevant existing
literature and background to be used in the sequel. Part
{\bf A} of Theorem \ref{main} is proved in Section
\ref{Sec:ProofGPE}.
Aspects of the (non) minimality of the realization in
\eqref{eq:L} are addressed in Section \ref{Sec:4}. On
the way, we introduce a test for the
non-minimality of a state space realization of an
arbitrary system (vanishing at infinity) by examining
common eigenvalues between the \mbox{$(n+p)\times(n+p)$}
system matrix $L$ and its $n\times n$ submatrix $A$, see
\eqref{StateSpace}. Part {\bf B} of Theorem \ref{main} is
proved in Section \ref{Sec:Odd}. In Section
\ref{Sec:StatFeedArbitrary} we relate minimal realization
of an arbitrary system  (vanishing at infinity) with the
ability of moving  its poles through static output feedback.
As a sample application of Theorem \ref{main}, in Section
\ref{Sec:StatFeedGPE} we study the effect of static output
feedback on $\mathcal{GPE}$ and ${\mathcal Odd}$ systems
and then prove Proposition \ref{FeedBackSpectFact}.

\section{Background and perspective}
\label{Sec:2}
\setcounter{equation}{0}

In this section we state known results to be used in the
sequel. Generalized positive functions were introduced in
the context of the Positive Real Lemma (PRL), see \cite{AM}
and references therein\begin{footnote}{The original
formulation was real. The case we address is in fact
{\em generalized} positive and {\em complex}, but we wish
to adhere to the commonly used term: Positive Real
Lemma.}\end{footnote}. Applications of $\mathcal{GP}$
functions to electrical networks appeared in \cite{IO},
and to control in \cite{LSC}, where they first casted in
a Linear Matrix Inequality (LMI) framework, see e.g.
\cite{BGFB} for more information on LMI. For more
application of the generalized PRL, see \cite{HSK}.
\vskip 0.2cm

We now recall in three characterizations of $\mathcal{GP}$
functions. We start with the Positive Real Lemma (PRL) as
first presented in \cite[Theorem 1]{DDGK} (up to substituting
the real setting by a complex one)\begin{footnote}{The
original formulation was real. The case we address is in
fact {\em generalized} positive and {\em complex}, but we
wish to adhere to the commonly used
term: Positive Real Lemma.}\end{footnote}. For sample of
other versions the PRL see e.g. \cite[Theorem 15.2]{BGKR},
\cite[Section 3A]{HSK} and the discussions in 
\cite[Section 2]{AL2} and in \cite[p. 34]{BGFB}.
\vskip 0.2cm

\begin{lemma}\label{gpLemma}
Let $L$ be a $(n+p)\times(n+p)$ realization matrix of a
$~p\times p$-valued rational function $\Psi(s)$
of McMillan degree $~q$, see \eqref{StateSpace}.

(I) $\Psi\in\mathcal{GP}$ if there exists a Hermitian
non-singular $~\hat{H}\in\C^{n\times n}$ so that
\begin{equation}\label{Lyap}
HL+L^*H\in\overline{\mathbb P}_{n+p}\quad\quad\quad
H={\rm diag}\{\hat{H},~I_p\}.
\end{equation}
If in addition the realization is minimal, i.e. $q=n$,
the converse is true as well: $~\Psi\in\mathcal{GP}$
implies \eqref{Lyap}.

(II) If in part (I) $~-\hat{H}\in{\mathbb P}_n$ then
$~\Psi\in{\mathcal P}$.
\end{lemma}
\vskip 0.2cm

It is of interest to recall that in
\cite[Section 7]{AL2} we pointed out that if $L$ satisfies
\eqref{Lyap}, whenever non-singular, also
$L^{-1}$ satisfies \eqref{Lyap}, with the same $H$.
\vskip 0.2cm

As a second characterizations of $\mathcal{GP}$
functions recall that a $p\times p$-valued function
$\Psi(s)$ belongs to $\mathcal{GP}$ if and only if
it can be factored as
\begin{equation}\label{eq:factor}
\Psi(s)=G(s)P(s)G^{\#}(s),
\end{equation}
where $~G,P$ are $p\times p$-valued, $G$ analytic in $\C_-$
and $P\in\mathcal P$. See for instance \cite{Ge} for the
corresponding result in the setting of functions meromorphic
in the open unit disk rather than in the right open
half-plane, see the discussion in
\mbox{\cite[Section 1]{AL1}}. To present the third
characterizations of $\mathcal{GP}$ functions we briefly
mention that a rational function $\Psi$ is in $\mathcal{GP}$
if and only if the kernel
\[
\frac{\Psi(s)+\Psi(w)^*}{s+w^*}
\]
has a finite number of negative squares in its domain of
definition in $\mathbb C_+$. This equivalent characterization
of (not necessarily rational) functions of the form
\eqref{eq:factor} appeared for the scalar case in \cite{DHdS}
and \cite{DLLS1} and extended in \cite{L1}, \cite{L2},
\cite{L3}. The significance of \eqref{eq:factor} to scalar
rational $\mathcal{GP}$ functions was recently treated in
\cite{AL1} and \cite{AL3}, where a more complete survey of the
literature can also be found.
\vskip 0.2cm

The following result is crucial to our construction.

\begin{theorem}\label{Pr:GPEchar}
Let $\Psi(s)$ be a matrix valued rational function.
The following are equivalent
\begin{itemize}
\item[(i)~~~]{}$\Psi\in\mathcal{GPE}$

\item[(ii)~~]{}$\Psi\in\mathcal{GP}$ and in addition
$\Psi(i\omega)=\Psi(i\omega)^*$ a.e. $\omega\in\R$.

\item[(iii)~]{}$\Psi(s)=G(s)G(s)^{\#}$ and without loss of
generality $G$ can be chosen so that $G$ and $G^{-1}$ are
analytic in $\C_-$.
\end{itemize}
\end{theorem}

{\bf Proof}\quad
The equivalence of (i) and (ii) follows from the definition
of $\mathcal{GPE}$ functions, see \eqref{DefGPE}.

The fact that $~(ii)~\Longrightarrow (iii)~$ is deep and was
established in \cite[Theorem 2]{Y}.

The fact that $~(iii)~\Longrightarrow (i)~$ is straightforward.
\qed
\vskip 0.2cm

The problem in item (iii) of finding $GG^{\#}$ out of
$\Psi\in\mathcal{GPE}$ is known as ~{\em spectral factorization}
whenever $\Psi$ and $\Psi^{-1}$ are analytic on the
imaginary axis and if this restriction is relaxed, ~{\em pseudo
spectral factorization}.
As sample references on spectral factorization, see e.g.
\cite[Chapter 9]{BGKR} \cite{Fu} \cite[Section 19.3]{LR} and
\cite{Ran}. On pseudo-spectral factorization, see e.g.
\cite[Chapter 10]{BGKR}, \cite{Roo} and the earlier work of
Youla on which we rely, \cite{Y}. The gap between spectral 
and pseudo spectral factorizations is addressed in 
\cite[Chapers 9, 10]{BGKR}. Here, Proposition
\ref{FeedBackSpectFact} is restated and proved in
Proposition \ref{FeedBackSpectFact1} below.
\vskip 0.2cm

For completeness we mention that
the spectral factorization problem has been extended to
the case where $F\in{\mathcal Even}$ can be factored to
\begin{center}
$F(s)=G(s)\cdot{\rm diag}\{I,~-I\}\cdot{G}^{\#}(s)$
\end{center}
see e.g. \cite{KLR}, \cite{LPR} and \cite[Part VII]{BGKR}.
\vskip 0.2cm

Note that $\Psi\in\mathcal{GPE}$ if and only if one
substitutes in \eqref{eq:factor} $\Psi=\Psi^{\#}$ thus
it is equivalent to having in \eqref{eq:factor}
$P(s)\equiv P$ for some~ {\em constant}~ positive
semidefinite $P$. This conforms well with item (iii)
in Theorem \ref{Pr:GPEchar}.
\vskip 0.2cm

The following result from \cite[Theorem 4.1]{AG1}, (see
also \cite[Proposition 2.1]{Ran}) is adapted to our
framework.
\begin{theorem}\label{Th:CharH}
Let $F(s)$ be a matrix valued rational function admitting
a state space realization $L$ as in \eqref{StateSpace} of
McMillan degree $q$. If
\begin{equation}\label{eq:CharH}
ML+(ML)^*=0_{n+p}
\end{equation}
where
\[
M:={\rm diag}\{H,~iI_p\},
\]
for some Hermitian non-singular $H\in\C^{n\times n}$, then
\begin{equation}\label{eq:Herm}
F(i\omega)=F(i\omega)^*\quad a.e. \quad\quad\omega\in\R.
\end{equation}
Conversely, let $L$ be a minimal realization of $F(s)$
in \eqref{eq:Herm}, i.e. $q=n$. Then $L$ satisfies
\eqref{eq:CharH}.
\end{theorem}

Note that $M$ is not Hermitian.
\vskip 0.2cm

Clearly, if $F$ satisfies \eqref{eq:Herm}, so does $-F$ and
also $F^{-1}$ (whenever the determinant is not
identically zero). It is less obvious that
some of this properties hold for the
system matrix $L$. Namely, if $L$ satisfies \eqref{eq:CharH}
with $M={\rm diag}\{H,~~iI_p\}$, then this holds for $-L$
and whenever exist, $\pm{L}^{-1}$ satisfies \eqref{eq:CharH}
with $M={\rm diag}\{-H,~~iI_p\}$.
\vskip 0.2cm

We conclude this section with recalling in a technical result
on the state space realization of a composition of a pair of
rational functions (series or cascade connection of systems
in Electrical Engineering terminology). Namely $F_{\alpha}(z)$,
$F_{\beta}(z)$ are of compatible dimensions and
$F_{\gamma}(z)$ is obtained by
\begin{equation}\label{eq:Fb}
F_{\gamma}(z)=F_{\alpha}(z)F_{\beta}(z).
\end{equation}
Assuming the state-space realization of each $F_{\alpha}(z)$
and $F_{\beta}(z)$ is known, one can construct a realization of
the resulting $F_{\gamma}(z)$, see e.g.
\cite[Subsection 8.3.3]{Ka}, \cite[Eq. (4.15)]{Wa}.

\begin{observation}\label{Ob:series}
Given $l\times q$ and $q\times r$ valued rational functions
$F_{\alpha}(z)$, $F_{\beta}(z)$, respectively, admitting
state space realization
\[
F_{\alpha}(z)=C_{\alpha}(zI-A_{\alpha})^{-1}B_{\alpha}
+D_{\alpha}\quad\quad\quad
A_{\alpha}\in\C^{p_{\alpha}\times p_{\alpha}}
\]
\[
\begin{matrix}
B_{\alpha}\in\C^{p_{\alpha}\times q}&~&~&
C_{\alpha}\in\C^{l\times p_{\alpha}}&~&~&
D_{\alpha}\in\C^{l\times q}.
\end{matrix}
\]
\[
F_{\beta}(z)=C_{\beta}(zI-A_{\beta})^{-1}B_{\beta}+D_{\beta}
\quad\quad\quad A_{\beta}\in\C^{p_{\beta}\times p_{\beta}}
\]
\[
\begin{matrix}
B_{\beta}\in\C^{p_{\beta}\times r}&~&~&
C_{\beta}\in\C^{q\times p_{\beta}}&~&~&
D_{\beta}\in\C^{q\times r}
\end{matrix}
\]
The system matrix $L_{\gamma}$ associated with the realization
of $F_{\gamma}(z)$ in \eqref{eq:Fb} is given by,
%$$
%L_{\gamma}=\left(\partition
%\matrix
%%\begin{smallmatrix}
%A_{\alpha}&B_{\alpha}C_{\beta}&\vdots&B_{\alpha}D_{\beta}\\
%0&~A_{\beta}&\vdots&B_{\beta}\\
%\cdots&\cdots&\cdots&\cdots\\
%C_{\alpha}&~D_{\alpha}C_{\beta}&\vdots&D_{\alpha}D_{\beta}
%\end{smallmatrix}
%\endmatrix
%\vdashed 3:5\hdashed 4:6\right)
%$$
\[
L_{\gamma}=\left(
\begin{array}{rr|r}
%{smallmatrix}
A_{\alpha}&B_{\alpha}C_{\beta}
%&\vdots
&B_{\alpha}D_{\beta}\\
0~~~~~&A_{\beta}
%&\vdots
&B_{\beta}\\ \hline
%\cdots&\cdots&\cdots&\cdots\\
C_{\alpha}&~D_{\alpha}C_{\beta}
%&\vdots
&D_{\alpha}D_{\beta}
\end
%{smallmatrix}
{array}\right).
\]
\end{observation}

\section{Proof of Theorem \ref{main} {\bf A}}
\label{Sec:ProofGPE}
\setcounter{equation}{0}

{\mbox{\boldmath $(ii)~\Longrightarrow (i)$}}

Substituting in Theorem \ref{Pr:GPEchar} $~(ii)~\Longrightarrow (i)$
a combination of Lemma \ref{gpLemma} (the PRL) and
Theorem \ref{Th:CharH}, establishes this part.
\vskip 0.2cm

{\mbox{\boldmath $(i)~\Longrightarrow (iii)$}}

First, substitute in \eqref{StateSpace} $F(s)=G(s)$
so that its (not necessarily minimal) realization is
\begin{equation}\label{RealizG}
L_g=\left(\begin{array}{r|r}
%{smallmatrix}
\hat{A}
%&\vdots
&\hat{B}\\ \hline
%\cdots&\cdots&\cdots\\
\hat{C}
%&\vdots
&\hat{D}\end
%{smallmatrix}
{array}\right).
\end{equation}
Then the
corresponding realization of $G(s)^{\#}$ is given by,
\[
\left(\begin
%{smallmatrix}
{array}{r|r}-\hat{A}^*
%&\vdots
&-\hat{C}^*\\ \hline
%\cdots&\cdots&\cdots\\
\hat{B}^*
%&\vdots
&\hat{D}^*\end
%{smallmatrix}
{array}\right).
\]
Now, substituting in Observation \ref{Ob:series} $F_{\alpha}=G$
and $F_{\beta}=G^{\#}$ along with
Proposition \ref{Pr:GPEchar} establishes the structure in
\eqref{eq:L}, i.e.
\begin{equation}\label{eq:L1}
L=\left(\begin
%{smallmatrix}
{array}{rr|r}
\hat{A}&\hat{B}\hat{B}^*
%&\vdots
&\hat{B}\hat{D}^*\\
0&-\hat{A}^*
%&\vdots
&-\hat{C}^*\\ \hline
%\cdots&\cdots&\cdots&\cdots\\
\hat{C}&\hat{D}\hat{B}^*
%&\vdots
&\hat{D}\hat{D}^*\end
%{smallmatrix}
{array}\right).
\end{equation}
Finally, we have shown above how to construct a realization of
$\mathcal{GPE}$ function of \mbox{degree$=2n$} where $n$ is
arbitrary, see \eqref{eq:L1}. Note now that if
$\Psi\in\mathcal{GPE}$ is realized by
\mbox{$L=\left(\begin
%{smallmatrix}
{array}{r|r}A
%&\vdots
&B\\ \hline
%\cdots&\cdots&\cdots\\
C
%&\vdots
&D\end
%{smallmatrix}
{array}\right)$} of the form in
\eqref{eq:L1}, then the same $\Psi$ can also realized by
$\tilde{L}=\left(\begin
%{smallmatrix}
{array}{rr|r}A&*
%&\vdots
&B\\ 0&*
%&\vdots
&0\\ \hline
%\cdots&\cdots&\cdots&\cdots\\
C&*
%&\vdots
&D\end
%{smallmatrix}
{array}\right)$, where $~*~$ means
an arbitrary block of suitable dimensions. It may be the case
that one can not transform $\tilde{L}$ by change of
coordinates to the form in \eqref{eq:L1} (in particular the
dimension of the $\tilde{A}$ part in $\tilde{L}$ may be odd).
However, as all minimal realizations are similar,
see e.g. \cite[8. 29]{Be}, \cite[Theorem 2.4-7]{Ka},
one of them is of the form \eqref{eq:L1}.
\vskip 0.2cm

{\mbox{\boldmath $(iii)~\Longrightarrow (ii)$}}

We find it convenient to introduce the following four
intermediate conditions:
\vskip 0.2cm

\begin{equation}\label{(ii)a}
{\rm Condition}~~~ (ii)~~~ {\rm in}~~ {\rm Theorem}~~~
\ref{main}~~~{\rm\bf A}~~~{\rm is}~~~ {\rm satisfied}~~~~
{\rm where}~~~~H_1~~~{\rm and} ~~~H_2~~ {\rm share}~~~
{\rm the}~~~{\rm same}~~{\rm inertia}.
\end{equation}

\begin{equation}\label{(ii)b}
{\rm Condition}~~~ (ii)~~~{\rm in}~~{\rm Theorem}~~~
\ref{main}~~~{\rm\bf A}~~~{\rm is}~~~{\rm satisfied}~~~
{\rm with}~~~H_2~~~{\rm unitarily}~~~{\rm similar}~~~
{\rm to}~~~H_1.
\end{equation}

\begin{equation}\label{(ii)c}
{\rm Condition}~~~(ii) ~~~{\rm in}~~~{\rm Theorem}~~~
\ref{main}~~{\rm\bf A}~~~{\rm is}~~~{\rm satisfied}~~~
{\rm with}~~~H_2~~~{\rm unitarily}~~~{\rm similar}~~~
{\rm to}~~~H_1,~~~{\rm both}~~~{\rm involutions.}
\end{equation}

\begin{equation}\label{(ii)d}
{\rm Condition} ~~(ii) ~~{\rm in}~~~{\rm Theorem}~~~
\ref{main} ~~{\rm\bf A} ~~{\rm is}~~~{\rm msatisfied}~~~~
{\rm with}\quad
H_1=\left(\begin{smallmatrix}0&I_n\\
I_n&0\end{smallmatrix}\right)\quad\quad
H_2=i\left(\begin{smallmatrix}0&-I_n\\
I_n&0\end{smallmatrix}\right).
\end{equation}
Trivially \eqref{(ii)d} $~\Longrightarrow~$ \eqref{(ii)c}
$~\Longrightarrow~$ \eqref{(ii)b} $~\Longrightarrow~$
\eqref{(ii)a} $~\Longrightarrow~$ (ii), so all we need
to show is that (iii) $~\Longrightarrow~$ \eqref{(ii)d}.
\vskip 0.2cm

Indeed, a straightforward calculation reveals that the
matrix $L$ in \eqref{eq:L1} satisfies the conditions in
\eqref{(ii)d} where \mbox{$HL+L^*H=
\left(\begin{smallmatrix}0\\ \sqrt{2}B\\ \sqrt{2}D
\end{smallmatrix}\right)
\left(\begin{smallmatrix}0\\ \sqrt{2}B\\ \sqrt{2}D
\end{smallmatrix}\right)^*$}. Verifying \eqref{eq:CharH}
is immediate so the proof is complete.
\qed
\vskip 0.2cm

Note that the matrix $L$ in \eqref{eq:L1} has a
special symmetry, i.e.
\begin{equation}\label{eq:symmetry}
\left(\begin{smallmatrix}
0&-I_n&0\\ I_n&0&0\\ 0&0&I_p\end{smallmatrix}\right)
\left(\begin{smallmatrix}\hat{A}&~\hat{B}\hat{B}^*
&\hat{B}\hat{D}^*\\ 0&-\hat{A}^* &-\hat{C}^*\\
\hat{C}&~\hat{D}\hat{B}^* &\hat{D}\hat{D}^*
\end{smallmatrix}\right)=\left(\begin{smallmatrix}
0&\hat{A}^*&\hat{C}^*\\
\hat{A}&\hat{B}\hat{B}^*&\hat{B}\hat{D}^*\\
\hat{C}&\hat{D}\hat{B}^*&\hat{D}\hat{D}^*
\end{smallmatrix}\right).
\end{equation}
In particular, the $A$ matrix in \eqref{eq:L} (the
upper left bock in $L$) has a Hamiltonian
structure\begin{footnote}{sometimes called
{\em $H$-skew-Hermitian}, see e.g.
\cite[sub-section 3.2]{LR}.}\end{footnote},
see also \eqref{eq:Acl} below. This will be
exploited in the proof Proposition
\ref{FeedBackSpectFact1} below.
\vskip 0.2cm

We now examine an aspect of Theorem \ref{main} {\bf A}:
Above we have {\em indirectly} proved that
\hfill{{\bf A} (ii) $\Longrightarrow$ \eqref{(ii)a}
$\Longrightarrow$ \eqref{(ii)b} $\Longrightarrow$
\eqref{(ii)c} $~\Longrightarrow~$ \eqref{(ii)d}
$\Longrightarrow$ {\bf A} (iii).}

We now illustrate the fact that directly showing these
implications is not obvious.
\vskip 0.2cm

\eqref{(ii)d} $~\Longrightarrow~$ (iii).

Let $\Delta_l=\left(\begin
%{smallmatrix}
{array}{rr|r}0~~~&P_2
%&\vdots
&0\\ P_1&0
%&\vdots
&0\\ \hline
%\cdots&\cdots&\cdots&\cdots\\
0~~~&0
%&\vdots
&P_3\end
%{smallmatrix}
{array}\right)$
where $P_1, P_2\in\overline{\mathbb P}_n$
and $P_3\in\overline{\mathbb P}_l$ are arbitrary.
Then $\Delta_l$ satisfies the conditions of
\eqref{(ii)d} with the same $H$ and $M$ as in \eqref{(ii)d}.
This in turn implies that also $L+\Delta_l$ satisfies
the same conditions. However, $L+\Delta_l$ is no longer
of the form of \eqref{eq:L}. This is illustrated by the
following example.
\vskip 0.2cm

Let $L_1=\left(\begin
%{smallmatrix}
{array}{rr|r}-1&1
%&\vdots
&~1\\ ~0&1
%&\vdots
&-1\\ \hline ~
%\cdots&\cdots&\cdots&\cdots\\
1&1
%&\vdots
&~1\end
%{smallmatrix}
{array}\right)$ of the of the form of \eqref{eq:L} be a
realization the rational $\mathcal{GPE}$ function
$\psi_1(s)=\frac{1}{1-s^2}+1$. Now in $\Delta_l$ substitute
$P_1=1$, $P_2=2$ and $P_3=3$, i.e.
$\Delta_l=\left(\begin
%{smallmatrix}
{array}{rr|r}0&2
%&\vdots
&0\\ 1&0
%&\vdots
&0\\ \hline
%\cdots&\cdots&\cdots&\cdots\\
0&0
%&\vdots
&3\end
%{smallmatrix}
{array}\right)$. Indeed, $\Delta_l$
satisfies condition \eqref{(ii)d} and so does
$L_2:=L_1+\Delta_l=\left(\begin
%{smallmatrix}
{array}{rr|r}-1&3
%&\vdots
&~1\\~ 1&1
%&\vdots
&-1\\ \hline ~
%\cdots&\cdots&\cdots&\cdots\\~
1&1
%&\vdots
&~4\end
%{smallmatrix}
{array}\right)$. Although $L_2$ is
not of the form of \eqref{eq:L}, it is a state space
realization of the rational $\mathcal{GPE}$ function
$\psi_2(s)=\frac{4}{4-s^2}+4.$ Indeed, $\psi_2(s)$ can also
be realized by $L_2'=\left(\begin
%{smallmatrix}
{array}{rr|r}-2~~~~~~~~&4~~~~~
%&\vdots
&4 \\~0~~~~~~~~&~~2
%&\vdots
&2-\sqrt{5}\\ \hline
%\cdots&\cdots&\cdots&\cdots\\
\sqrt{5}-2&4~~~~~
%&\vdots
&4\end
%{smallmatrix}
{array}\right)$ which is of
the form of \eqref{eq:L}.
\vskip 0.2cm

\eqref{(ii)b} $~\Longrightarrow~$ \eqref{(ii)c}.

Consider the condition for the submatrix that
$H_1A+A^*H_1\in\overline{\mathbb P}_n~$ and
\mbox{$H_2A+A^*H_2=0_n$.} Take $A=\left(\begin{smallmatrix}
0&5\\ 1&0\end{smallmatrix}\right)$. $H_1={\rm diag}\{1,~-5\}$
and $H_2=\left(\begin{smallmatrix}0&\sqrt{5}\\
\sqrt{5}&-4\end{smallmatrix}\right)$. Then the conditions
are satisfied with $H_1$ and $H_2$ unitarily similar, but neither
is an involution.
\vskip 0.2cm

\eqref{(ii)a} $~\Longrightarrow~$ \eqref{(ii)b}.

Consider the condition for the submatrix that
$H_1A+A^*H_1\in\overline{\mathbb P}_n~$ and
\mbox{$H_2A+A^*H_2=0_n$.} Take $A=\left(\begin{smallmatrix}
0&5\\ 1&0\end{smallmatrix}\right)$. $H_1={\rm diag}\{1,~-5\}$
and $H_2=\left(\begin{smallmatrix}1&~~1\\
1&-3\end{smallmatrix}\right)$. Then the conditions
are satisfied with $H_1$ and $H_2$ sharing the same inertia,
but they are not unitarily similar.

\section{minimality of the realization}
\label{Sec:4}
\setcounter{equation}{0}

We now address the question of minimality of the obtained
realization. We first resort to the following observation
which goes beyond the scope of this work.

\begin{lemma}\label{LemmaMinReal}
Let $L\in\C^{(n+p)\times (n+p)}$ be realization of a
$p\times p$-valued rational function $F(s)$ of McMillan
degree $n$. Namely,
\[
\begin{matrix}
F(s)=C(sI-A)^{-1}B+D&~&~&
L:=\left(\begin{smallmatrix}A&~B\\C&~D
\end{smallmatrix}\right)\end{matrix},
\]
see \eqref{StateSpace}. If the realization is not minimal,
there exists $\lambda\in{\rm spect}(A)$ so that
$\lambda\in{\rm spect}(L)$ for all $D$.
\vskip 0.2cm

If ${\rm rank}(B)={\rm rank}(C)$
or the number of Jordan blocks in $A$ is less
or equal to $\min({\rm rank}(B),~{\rm rank}(C))$
then the converse is true as well.
\end{lemma}

{\bf Proof}\quad
We first show that if the realization is not minimal
there exists \mbox{$\lambda\in{\rm spect}(A)$} so that
$\lambda\in{\rm spect}(L)$ for all $D$.
\vskip 0.2cm

If the realization $L$ is not controllable, there exists
\mbox{$0\not=\hat{v}\in\C^n$} so that
\mbox{$\hat{v}^*(A-\lambda{I}_n)=0$,} for some $\lambda\in\C$
and in addition \mbox{$\hat{v}^*B=0$}, see
e.g. \mbox{\cite[Theorem 3.3.1]{AV}},
\mbox{\cite[Theorem 2.4-8]{Ka}},
\mbox{\cite[Theorem 4.3.3]{LR}}. This implies that
$(n+p)$-dimensional vector
$v:=${\mbox{\tiny$\begin{pmatrix}\hat{v}\\ 0\end{pmatrix}$}}
satisfies, $v^*(L-\lambda{I_{n+p}})=0$ with the same $\lambda$.
Namely, $\lambda\in{\rm spect}(A)\bigcap{\rm spect}(L)$.
\vskip 0.2cm

Similarly, if the realization $L$ is not observable, there
exists $0\not=\hat{v}\in\C^n$ so that
$(A-\lambda{I}_n)\hat{v}=0$, for some $\lambda\in\C$ and
in addition $C\hat{v}=0$, see
e.g. \mbox{\cite[Theorem 3.3.6]{AV}},
\mbox{\cite[Theorem 2.4-8]{Ka}}.
This implies that $(n+p)$-dimensional vector
$v:=${\mbox{\tiny$\begin{pmatrix}\hat{v}\\ 0\end{pmatrix}$}}
satisfies, $(L-\lambda{I_{n+p}})v=0$ with the same $\lambda$.
\vskip 0.2cm

For the converse direction we now show that if the
realization $L$ is minimal one can construct $D$ so that
the matrix
\[
L-\lambda{I}_{n+p}=
\left(\begin{smallmatrix}A-\lambda{I}_n&~B\\
C&~D-\lambda{I}_p\end{smallmatrix}\right)
\]
is nonsingular for each \mbox{$\lambda\in{\rm spect}(A)$},
i.e ${\rm spect}(A)\bigcap{\rm spect}(L)=\emptyset$. To
this end, take $D$ so that
${\rm spect}(A)\bigcap{\rm spect}(D)=\emptyset$. Hence,
one can write
\[
L-\lambda{I}_{n+p}=\left(\begin{smallmatrix}I_n~~&
B(D-\lambda{I}_p)^{-1}\\0&I_p\end{smallmatrix}\right)
\left(\begin{smallmatrix}A-\lambda{I_n}+B(\lambda{I}_p
-D)^{-1}C&0\\
C&~~~D-\lambda{I}_p\end{smallmatrix}\right).
\]
In this case $\lambda\in{\rm spect}(L)$ if and only if
the matrix
\[
A-\lambda{I}_n+BKC\quad\quad{\rm with}
\quad\quad K:=(D-\lambda{I}_p)^{-1}
\]
is singular. The construction for a minimal realization
of such $K\in\C^{p\times p}$ so that $A-\lambda{I}_n+BKC$
is nonsingular for all \mbox{$\lambda\in{\rm spect}(A)$}
is established in Proposition \ref{Pr:OutputFeedback}
below, so the proof is complete.
\qed
\vskip 0.2cm

Admittedly, we do not know whether or not the conditions
on ${\rm rank}(B)$ and ${\rm rank}(C)$ in the above
result are inherent to the problem or just a by-product
of the technique we employed.
\vskip 0.2cm

We now illustrate the fact that it may be the case that
\mbox{${\rm spect}(A)\bigcap{\rm spect}(L)\not=\emptyset$}
for some $D$, although a realization is minimal. Indeed,
take $A=0$ (scalar) \mbox{$B^*, C\in\C^2$,} both non-zero,
and \mbox{$D\in\C^{2\times 2}$}. Namely the following
$3\times 3$
system matrix \mbox{$L=\left(\begin{smallmatrix}0&B\\
C&D\end{smallmatrix}\right)$} is a minimal (degree=1)
realization of a $2\times 2$-valued rational function
\mbox{$F(s)=\frac{1}{s}CB+D$.}
Obviously, minimality is independent of $D$. However, if
for example \mbox{$F(s)=(\frac{1}{s}+\gamma)CB$,} namely
$D={\gamma}CB$, where $\gamma\in\C$ is arbitrary, both $A$
and $L$ are singular, i.e. have a common eigenvalue (=zero).
\vskip 0.2cm

The special structure of $D$ in the above example suggests
that if a realization is minimal, the matrices $L$ and $A$
typically do not share common eigenvalues.
\vskip 0.2cm

One can now exploit the special structure of the realization
$\mathcal{GPE}$ functions in \eqref{eq:L1},
\eqref{eq:symmetry}, in applying 
Lemma \ref{LemmaMinReal} to obtain the following.

\begin{observation}\label{spectL-A}
Let $L_g(s)$ in \eqref{RealizG} be a realization of
$p\times p$-valued rational function $G(s)$ and without
loss of generality assume that
${\rm spect}(\hat{A})\subset\overline{\C_+}$.

Let $\tilde{G}(s)=G(s)-\hat{D}$, i.e.
$\tilde{G}(s)=\hat{C}(sI_n-\hat{A})^{-1}\hat{B}$.
The following is true.
\vskip 0.2cm

{\bf A}.  Minimality of the realizations of $G(s)$ and
of $\tilde{G}(s)$ is equivalent.
\vskip 0.2cm

{\bf B}. The following are equivalent.
\begin{itemize}

\item[(i)~~~]{}The realization of $\tilde{G}(s)$ is not
minimal.

\item[(ii)~~]{}There exists $\lambda\in{\rm spect}(\hat{A})$
so that for all $\hat{D}$,
$\lambda\in{\rm spect}(L_g)$ in \eqref{RealizG}.

\item[(iii)~]{}For all $\hat{D}$ the realization $L$ in
\eqref{eq:L1}, of $G(s)G(s)^{\#}$ is not minimal.

\item[(iv)~]{}There exists \mbox{$\lambda\in{\rm spect}(
\hat{A})\bigcup{\rm spect}(-\hat{A}^*)$} so that for all
$\hat{D}$,  $\lambda\in{\rm spect}(L)$ in \eqref{eq:L1}.
\end{itemize}
\vskip 0.2cm

{\bf C}. Assume that the realization of $\tilde{G}(s)$ is
minimal. If ${\rm spect}(\hat{A})\not\subset{i}\R$, one
can always find $\hat{D}$ so that the realization of
$G(s)G(s)^{\#}$ in \eqref{eq:L1}, is not minimal.

\end{observation}

{\bf Proof}\quad {\bf A}. As minimality of realization is
equivalent
to controllability and observability, the claim is obvious.
\vskip 0.2cm

{\bf B}. First note 
that if $C(sI-A)^{-1}B+D$ is a realization of a
$\mathcal{GPE}$ function of the form \eqref{eq:L1}, then
from \eqref{eq:symmetry} it in particular follows that
$C^*=${\mbox{\tiny$\begin{pmatrix}0&-I_n\\ I_n&~~0
\end{pmatrix}$}}$B$, which implies
${\rm rank}(B)={\rm rank}(C)$. Thus, one can
apply Lemma \ref{LemmaMinReal} which now implies
\mbox{(i) $\Longleftrightarrow$ (ii)} ~and~
\mbox{(iii) $\Longleftrightarrow$ (iv).}
\vskip 0.2cm

{\bf (i) $\Longrightarrow$ (iii)}

If the realization of $\tilde{G}(s)$ is not minimal, assume
it is not observable. Namely, there exists
\mbox{$0\not=v\in\C^n$} so that {\mbox{\tiny$
\begin{pmatrix}v\\ 0\end{pmatrix}$}} is an
$(n+p)$-dimensional vector of $L_G$ in \eqref{RealizG}
(corresponding to an eigenvalue $\lambda$). Then,
independent of $\hat{D}$,
{\mbox{\tiny$\begin{pmatrix}v\\ 0\end{pmatrix}$}} is an
$(2n+p)$-dimensional vector of $L$ in \eqref{eq:L1}
(corresponding to an eigenvalue $\lambda$). Thus, this
realization is not observable. 
(controllability can be similarly treated).
\vskip 0.2cm

{\bf (iii) $\Longrightarrow$ (i)}

Assume that the realization $L$ in
\eqref{eq:L1} is not minimal for all $\hat{D}$. If it is not
observable (controllability can be similarly treated), there
exist $v_1,v_2\in\C^{n}$ so that
{\mbox{\tiny$\begin{pmatrix}v_1\\ v_2\\ 0\end{pmatrix}$}}
is a non-zero $(2n+p)$-dimensional eigenvector of $L$,
corresponding
to an eigenvalues $\lambda$. As this holds for all $\hat{D}$,
it implies that $\hat{B}^*v_2=0$. Now if $v_2\not=0$ then it is
also a left eigenvector of $\hat{A}$ (corresponding to an
eigenvalues $-\lambda^*$) so the pair $\hat{A}, \hat{B}$ is
not controllable. If $v_2=0$ then $0\not=v_1$ is an eigenvector
of $\hat{A}$ (corresponding to an eigenvalues $\lambda$) and
$\hat{C}v_1=0$, thus the pair $\hat{A}, \hat{C}$ is not
observable.
\vskip 0.2cm

{\bf C}. Let $\lambda\in\C_+$ be an eigenvalue of $\hat{A}$
and let \mbox{$0\not=v\in\C^n$} be a corresponding
left eigenvector. Namely, $v^*(\hat{A}-\lambda{I})=0$,
i.e. $(\hat{A}^*-{\lambda}^*I)v=0$. Minimality of
the realization of $\tilde{G}(s)$ implies that
$v^*\hat{B}\not =0$ and $\hat{C}v\not=0$.
\vskip 0.2cm

By construction the matrix $(\hat{A}+\lambda^*I)$ is
nonsingular, so the non-zero $(2n+p)$-dimensional vector
{\mbox{\tiny$\begin{pmatrix}
-(\hat{A}+\lambda^*I)^{-1}\hat{B}\hat{B}^*v\\
v\\ 0\end{pmatrix}$}} is well defined. Now,
for $\hat{D}=\hat{C}(\hat{A}+\lambda^*I)^{-1}\hat{B}$
this is an eigenvector of $L$ in \eqref{eq:L1}
(corresponding to an eigenvalue $-\lambda^*)$ and thus
this realization is not observable
(controllability can be similarly treated).
\qed
\vskip 0.2cm

We conclude this section by illustrating the fact that the
condition ${\rm spect}(A)\not\subset{i}\R$ in item {\bf C}
is essential.

\begin{example}\label{Ex:GPEImagSpect}
{\rm
If ${\rm spect}(\hat{A})\subset{i\R}$, it may
be that the realization $L$ in \eqref{eq:L1} is minimal for
all $D$. Take for example the scalar function
$\Psi_1(s)=\frac{1}{s}$ which is positive and odd.
Its minimal state space realization is
$\Psi_1(s)=\hat{C}(sI-\hat{A})^{-1}\hat{B}$ with
$\hat{A}=0$, $\hat{B}=1$ and $\hat{C}=1$.
\vskip 0.2cm

Let $\Psi_2\in\mathcal{GPE}$ be of the form
\mbox{$\Psi_2(s)=\Psi_1(s)\Psi_1^{\#}(s)=-\frac{1}{s^2}$}.
Its realization is of the form of \eqref{eq:L},
given by
\begin{equation}\label{eq:ExImagSpect}
A=\left(\begin{smallmatrix}0&1\\0&0\end{smallmatrix}\right)
\quad\quad\quad B=\left(\begin{smallmatrix}~~0\\
-1\end{smallmatrix}\right)\quad\quad\quad C=(1,~~0)
\quad\quad\quad D=0.
\end{equation}
In particular, this realization is minimal.
\vskip 0.2cm

Consider now the $\mathcal{GPE}$ function $\Psi_2(s)+D$,
where $D\geq 0$ is a parameter. For all $D$, the realization
is minimal, i.e. of degree 2.
}
\qed
\end{example}

\section{Odd systems}
% - Proof of Theorem \ref{main} {\bf B}}
\label{Sec:Odd}
\setcounter{equation}{0}

{\bf Proof of Theorem \ref{main} {\bf B}}

{\bf (ii) $\Longrightarrow $ (iii)}

Let $\hat{H}$ has $\nu$ eigenvalues in $\C_-$ and $n-\nu$
eigenvalues in $\C_+$. Then up to similarity on $L$ (and
$A$) and congruence on $H$ (and $\hat{H}$) one can
take $\hat{H}={\rm diag}\{-I_\nu~,~I_{n-\nu}\}$ (and thus
$H={\rm diag}\{-I_{\nu}, I_{n-\nu},~I_p\}$). Substituting
in $HL+L^*H=0_{n+p}$ yields \eqref{RealOdd}.
\vskip 0.2cm

{\bf (iii) $\Longrightarrow $ (i)}

A straightforward calculation yields $\Psi(s)=C(sI_n-A)^{-1}B+D$.
We shall concentrate on $\Psi(s)_{|_{s=i\omega}}$ for $\omega\in\R$.
We first conformally partition
\[
(i\omega{I_n}-A)^{-1}=\left(\begin{smallmatrix}
i\omega{I_{\nu}}-T_1&-\hat{A}\\-\hat{A}^*&
i\omega{I_{n-\nu}}-T_2\end{smallmatrix}\right)^{-1}:=
\left(\begin{smallmatrix}W(i\omega)&X(i\omega)\\
Y(i\omega)&Z(i\omega)\end{smallmatrix}\right).
\]
We thus formally have,
\begin{equation}\label{rational}
\Psi(i\omega)=B_1^*W(i\omega)B_1+B_1^*X(i\omega)B_2
-B_2^*Y(i\omega)B_1-B_2^*Z(i\omega)B_2+T_3.
\end{equation}
As for all $\omega\in\R$ the matrix
$~i\cdot{\rm diag}\{I_{\nu}~,~-I_{n-\nu}\}(i\omega{I_n}-A)$
is Hermitian, so is its inverse
\[
(i\omega{I_n}-A)^{-1}\cdot i\cdot{\rm diag}\{-I_{\nu}~,~
I_{n-\nu}\}=\left(\begin{smallmatrix}-iW(i\omega)&iX(i\omega)\\
-iY(i\omega)&iZ(i\omega)\end{smallmatrix}\right).
\]
Namely, \mbox{$W(i\omega)=-W(i\omega)^*:=T_4(i\omega)$},
\mbox{$Z(i\omega)=-Z(i\omega)^*:=T_5(i\omega)$}
and \mbox{$X(i\omega)=Y(i\omega)^*$.}
Substituting now in \eqref{rational} yields
\[
\Psi(i\omega)=B_1^*T_4(i\omega)B_1+B_1^*X(i\omega)B_2-
(B_1^*X(i\omega)B_2)^*-B_2^*T_5(i\omega)B_2+T_3.
\]
Namely \eqref{eq:SkewHerm} is obtained and this part of
the claim is established.
\vskip 0.2cm

{\bf (i) $\Longrightarrow $ (ii)}

First consider the condition in (ii) and note that
if one denotes $N:={\rm diag}\{I_n,~iI_p\}$ then
\begin{equation}\label{PRLodd}
0=HL+L^*H=HNN^*L+L^*NN^*H=M\hat{L}+\hat{L}^*M^*=
M\hat{L}+(M\hat{L})^*,
\end{equation}
where \mbox{$M:={\rm diag}\{\hat{H},~iI_p\}$} and
$\hat{L}=N^*L=${\mbox{\tiny$\begin{pmatrix}~~A&~~B\\
-iC&-iD\end{pmatrix}$}} i.e. if $L$ realizes

\mbox{$\Psi(s)=C(sI_n-A)^{-1}B+D$,} $\hat{L}$ realizes
$F(s):=-i\Psi(s)$.
It now follows that on the
imaginary axis $F(s)$ is Hermitian, see \eqref{eq:Herm}
and $\Psi(s)$ is skew-Hermitian, i.e.
$\Psi\in{\mathcal Odd}$.
\vskip 0.2cm

Assuming the realization is minimal, i.e. $q=n$, both
directions of Theorem \ref{Th:CharH} holds, so in the
proof is complete.
\qed
\vskip 0.2cm

Note that as ${\mathcal Odd}\subset\mathcal{GP}$,
Lemma \ref{gpLemma} (the PRL) is satisfied: Indeed
\eqref{PRLodd}
%\mbox{$HL+L^*H=0$} 
is a special case of \eqref{Lyap}.
\vskip 0.2cm

Of a special interest is the set $\mathcal{PO}$
of positive-odd functions (a.k.a Lossless or Foster), i.e.
${\mathcal Odd}\bigcap\mathcal{P}$. In electrical
networks theory they are associated with ~\mbox{L-C}~
circuits. For more details see e.g.
\cite[2.20, 2.36, 2.39, 7.33, 8.35, 8.36, 8.37]{Be},
\cite[pp. 12, 221, Theorem 2.7.4]{AV},
\cite[subsections 4.2, 5.1]{CL4}.
Recall also from \eqref{eq:factor} it follows
that $\Psi\in{\mathcal Odd}$ can always be factored as
$\Psi=GPG^{\#}$ with $P\in\mathcal{PO}$. Thus in a
sense, $\mathcal{PO}$ functions generate all ${\mathcal Odd}$
functions. The following is known, see e.g.
\cite[Eq. (5.2.6)]{AV}. It is immediate from part {\bf B}
of Theorem \ref{main} upon substituting in \eqref{RealOdd}
$\nu=n.$

\begin{corollary}
Let $\Psi(s)$ be a $p\times p$-valued rational function
so that $\lim\limits_{s\rightarrow\infty}\Psi(s)$ exists.
Let $L\in\C^{(n+p)\times(n+p)}$ be the associated system
matrix. If
\[
L=\left(\begin
%{smallmatrix}
{array}{r|r}T_n&B\\ \hline
B^*&T_p\end
%{smallmatrix}
{array}\right)
\quad\quad\quad\quad
\begin{smallmatrix}T_n\in\C^{n\times n}&~&T_n=-T_n^*\\
T_p\in\C^{p\times p}&~&T_p=-T_p^*\end{smallmatrix}
\]
Then, $\Psi\in\mathcal{PO}$.

If the realization is minimal, up to similarity,
the converse is true as well.
\end{corollary}

We conclude this section with the following
observation.
Clearly if $\Psi_1(s), \Psi_2(s)$ are $\mathcal{PO}$
functions of the same dimensions, then
$(\Psi_1-\Psi_2)\in{\mathcal Odd}$. Comparison between
the above Corollary and part {\bf B} of Theorem
\ref{main}, reveals that not every ${\mathcal Odd}$
function can be written as a difference of
$\mathcal{PO}$ of functions.

\section{Static output feedback - arbitrary functions}
\label{Sec:StatFeedArbitrary}
\setcounter{equation}{0}

Recall that applying a static output feedback to an
input-output system \mbox{$y(s)=F(s)u(s)$} ($u$ is
$m$-dimensional input and $y$ is $p$-dimensional output)
means taking $u=Ky+u'$ with $u'$ an auxiliary input and
$K$ a $m\times p$ {\em constant} matrix. The resulting
closed loop system is $y(s)=F_{\rm c.l.}(s)u'(s)$ with
$F_{\rm c.l.}=(I_p-FK)^{-1}F$.
\vskip 0.2cm

For simplicity, we adopt the common assumption that
$\lim\limits_{s\rightarrow\infty}F(s)=0$ (``strictly proper"
in engineering jargon). Thus, in \eqref{StateSpace} the
$(n+p)\times(n+m)$ realization matrix $L$ is of the form
$L=${\mbox{\tiny$\left(\begin
%{pmatrix}
{array}{r|r}A&B\\ \hline
C&0\end
%{pmatrix}
{array}\right)$}}. After applying a static output
feedback, the closed loop realization matrix $L_{\rm cl}$ is
\begin{equation}\label{Lcl}
L_{\rm c.l.}=\left(\begin
%{smallmatrix}
{array}{r|r}A_{\rm cl}&~~B\\ \hline
C&~~0\end
%{smallmatrix}
{array}\right)\quad\quad\quad
A_{\rm cl}:=A+BKC.
\end{equation}
The simplicity of the static output feedback has made it
very attractive. However, exploring its properties turned
out to be challenging. The most common associated problem
has been stabilization, namely guaranteeing (in the
continuous time case) that
${\rm spect}(A_{\rm cl})\subset\C_-$. This is illustrated
in a basic way in \cite[Section 3.1]{Ka} and for a sample
of more recent references see e.g. \cite{ElGOAiR},
\cite{HL}, \cite{KDS}, \cite{PFA}, \cite{SADG}. One can
go beyond stability, and in \cite[Proposition 8.1]{AL2}
we characterized systems which may turned to be
$\mathcal{GP}$, and in particular
$\mathcal{P}$, through static output feedback.
\vskip 0.2cm

Our first result here goes beyond the scope of this work.
To this end, we need the following notation.  Let $L$ be
realization of a $p\times m$-valued rational function
$F(s)$ of McMillan degree $n$ with
$\lim\limits_{s\rightarrow\infty}F(s)=0$. Namely,
\begin{equation}\label{RealF}
\begin{matrix}
F(s)=C(sI-A)^{-1}B&~&~&
L:=\left(\begin
%{smallmatrix}
{array}{r|r}A&~B\\ \hline C&~0
\end
%{smallmatrix}
{array}\right)\end{matrix},
\end{equation}
see \eqref{StateSpace}. I shall find it convenient to
denote $\beta:={\rm rank}(B)$, $\gamma:={\rm rank}(C)$.
Namely, there exist nonsingular matrices 
$R_b\in\C^{m\times m}$  and $R_c\in\C^{p\times p}$ so that
\mbox{$R_cF(s)R_b=${\mbox{\tiny$\begin{pmatrix}\hat{F}(s)&0\\
0&0_{(p-\gamma)\times(m-\beta)}\end{pmatrix}$}}} where
$\hat{F}(s)$ is $\gamma\times\beta$-valued.

We shall also
denote by $r$ the number of Jordan blocks in $A$.

\begin{proposition}\label{Pr:OutputFeedback}
Let $L$ be realization of $F(s)$ as in \eqref{RealF}.
If this realization is not minimal, there exists
\mbox{$\lambda\in{\rm spect}(A)$} so that for all
$K\in\C^{m\times p}$,
\mbox{$\lambda\in{\rm spect}(A_{\rm cl})$}, see
\eqref{Lcl}.
\vskip 0.2cm

If $\min(\beta,~\gamma)\geq r$
or $\beta=\gamma$, then the converse is true as well.
\end{proposition}

{\bf Proof}\quad
We first show that if the realization is not minimal then
there exists \mbox{$\lambda\in{\rm spect}(A)$} so that for
all $K$, \mbox{$\lambda\in{\rm spect}(A_{\rm cl})$.} If the
realization is not controllable, there exists a left
eigenvector $0\not=v_L\in\C^n$ (the subscript stands for
``left") so that
$v_L^*(A-\lambda{I}_n)=0$, for some $\lambda\in\C$ and in
addition \mbox{$v_L^*B=0$}, see
e.g. \mbox{\cite[Theorem 2.4-8]{Ka}},
\mbox{\cite[Theorem 4.3.3]{LR}}. This implies that
\[
v_L^*(A_{\rm cl}-\lambda{I}_n)=v_L^*(A+BKC-\lambda{I}_n)=
v_L^*(A-\lambda{I}_n)=0.
\]
Similarly, if the realization is not observable, there
exists $0\not=v_R\in\C^n$ (the subscript stands for ``right")
so that \mbox{$(A-\lambda{I}_n)v_R=0$,} for some
$\lambda\in\C$ and in addition \mbox{$Cv_R=0$}, see
e.g. \mbox{\cite[Theorem 2.4-8]{Ka}}.
This implies that
\[
(A_{\rm cl}-\lambda{I}_n)v_R=(A+BKC-\lambda{I}_n)v_R=
(A-\lambda{I}_n)v_R=0.
\]
\vskip 0.2cm

For the converse direction, assuming the realization is minimal,
we now construct $K\in\C^{m\times p}$ so that
\mbox{${\rm spect}(A)\bigcap{\rm spect}(A_{\rm cl})=\emptyset$.}
To this end, $K$ should be so that
\[
\begin{smallmatrix}
0\not=v_L\in\C^n&v_L^*(A-\lambda{I}_n)=0&\Longrightarrow&
v_L^*(A_{\rm cl}-\lambda{I}_n)=v_L^*BKC\not=0\\~\\
0\not=v_R\in\C^n&(A-\lambda{I}_n)v_R=0&\Longrightarrow&
(A_{\rm cl}-\lambda{I}_n)v_R=BKCv_R
\not=0
\end{smallmatrix}
\]
where the subscripts $L$ and $R$ stand for ``left" and
``right".
We shall find it convenient to denote \mbox{$K:=\delta\hat{K}$}
where \mbox{$0\not=\delta\in\C$} will be later determined.
Thus, we actually look for $\hat{K}$ (which will turn to
be an isometry if $m\geq p$, else coisometry) so that
\begin{equation}\label{CondK}
\begin{smallmatrix}
0\not=v_L\in\C^n&v_L^*(A-\lambda{I}_n)=0&\Longrightarrow&
v_L^* {B}\hat{K}C\not=0\\~\\0\not=v_R\in\C^n&(A-\lambda{
I}_n)v_R=0&\Longrightarrow&{B}\hat{K}Cv_R\not=0.
\end{smallmatrix}
\end{equation}
Consider now the singular values decomposition, see e.g.
\cite[Theorem 7.3.5]{HJ1}, of the matrices $B$ and $C$
\[
%\begin{matrix}
B=U_b\Sigma_bW_b
\quad \quad \quad \quad
C=U_c\Sigma_cW_c
\]
with $U_b,W_c\in\C^{n\times n}$, $W_b\in\C^{m\times m}$
and $U_c\in\C^{p\times p}$ all unitary and
\[
\Sigma_b=\left(\begin{smallmatrix}\hat{\Sigma}_b&0\\
0&~~~~~~~~~~~~~~0_{(n-\beta)\times(m-\beta)}
\end{smallmatrix}\right)
\quad\quad\quad\quad
\Sigma_c=\left(\begin{smallmatrix}\hat{\Sigma}_c
&0~~~~~~~~~~~~~~\\
0&~~~0_{(p-\gamma)\times(n-\gamma)}
\end{smallmatrix}\right)
\]
with $\hat{\Sigma_b}$ and $\hat{\Sigma_c}$ positive
diagonal of dimensions $\beta\times\beta$
and $\gamma\times\gamma$ respectively.
\vskip 0.2cm

If $\beta=\gamma$ then it is sufficient to take
\[
\hat{K}=W_b^*\left(\begin{smallmatrix}
I_{\beta}&0\\ 0&0_{(m-\beta)\times(p-\beta)}
\end{smallmatrix}\right)U_c^*
\]
so taht
\[
BKC={\delta}U_b\left(\begin{smallmatrix}
\hat{\Sigma}_b\hat{\Sigma}_c&0\\ 0&
0_{(n-\beta)\times(n-\beta)}
\end{smallmatrix}\right)W_b~.
\]
Minimality guarantees that \eqref{CondK} is satisfied.
\vskip 0.2cm

Next assuming that $\min(\beta,~\gamma)\geq r$ let
\[
\mathcal{V}_L:=\{B^*v_L~:~v_L^*(A-\lambda{I}_n)=0\}
\quad\quad\quad\quad
\mathcal{V}_R:=\{Cv_R~:~(A-\lambda{I}_n)v_R=0\}.
\]
Namely, \mbox{$\mathcal{V}_L\subset\C^m$} and
\mbox{$\mathcal{V}_R\subset\C^p$} can be spanned by
$r$ vectors and let \mbox{$K_L\in\C^{m\times r}$}
and \mbox{$K_R\in\C^{p\times r}$} be isometries
whose range is equal to the span of
\mbox{$\mathcal{V}_L$ and $\mathcal{V}_R$}, respectively.
Thus, if one takes
\[
\hat{K}=K_LK_R^*~,
\]
minimality guarantees that \eqref{CondK} is satisfied.
Thus, there is only the value of the scalar $\delta$ yet
to consider.
\vskip 0.2cm

If $A=\lambda{I_n}$ for some $\lambda\in\C$, minimality of
the realization implies that \mbox{$m\geq\beta= n$} and
\mbox{$p\geq\gamma=n$}. Thus, $\hat{K}=W_b^*U_c^*$ so any
non-zero $\delta$ will do. Hence, assume hereafter
\mbox{$A\not=\lambda{I_n}$.}
\vskip 0.2cm

In fact, we claim that for any $|\delta|$ sufficiently small,
$K$ satisfies the requirements. To this end
recall, see e.g. \cite[Corollary 7.3.8]{HJ1},
that if $M$ and $\Delta$ are two matrices of the same
dimensions then for all $k$,
\[
\|\Delta\|_2\geq |\sigma_k(M+\Delta)-\sigma_k(M)|
\]
where $\sigma_k(\cdot)$ are the respective singular
values $(\sigma_1\geq\sigma_2\geq~\cdots)$. For $j=1, 2,~\ldots$
let us denote by $\eta_j$ the smallest positive singular
value of \mbox{$A-\lambda_j{I}_n$} where
\mbox{$\lambda_j\in{\rm spect}(A)$}. Next let
$\eta:=\min(\eta_1~,~\eta_2~,~\ldots~)$. Taking,
\[
\delta=(2\eta\| B\|_2\| C\|_2)^{-1}
\]
guarantees that none of the positive singular values of
\mbox{$A-\lambda_j{I_n}$} with $j=1,~2,~\ldots$ was moved
to zero. Namely, $A_{\rm cl}-\lambda{I}_n$ is nonsingular
for all $\lambda\in{\rm spect}(A)$, so the construction
is complete.
\qed
\vskip 0.2cm

As already mentioned, we do not know whether or not the
above conditions on $\beta, \gamma$ are inherent to the
problem or just a by-product of the technique we employed.

\section{Static output feedback - $\mathcal{GPE}$ functions}
\label{Sec:StatFeedGPE}
\setcounter{equation}{0}

We start with the following question: Under what conditions
the $\mathcal{GP}$  class, and its subsets of $\mathcal{GPE}$
and ${\mathcal Odd}$,
are invariant under static output feedback.

\begin{proposition}\label{GPstaticOutputFeedback}
Consider the system $y(s)=\Psi(s)u(s)$ where
$\Psi\in\mathcal{GP}$ is $p\times p$-valued rational
function s.t. $\lim\limits_{s\rightarrow\infty}\Psi(s)=0$.
Consider the static output feedback \mbox{$u=Ky+u'$} where
$K$ is a constant $p\times p$ matrix with $u'$ auxiliary
input.
\vskip 0.2cm

{\bf A}. for all $K$ so that $-(K+K^*)\in\overline{\mathbb P}_p$
also the resulting closed loop system is in $\mathcal{GP}$.
\vskip 0.2cm

{\bf B}. If in addition $\Psi$ is even, i.e.
$\psi\in\mathcal{GPE}$ then, for all $K$ so that
$-K\in\overline{\mathbb P}_p$ also the resulting closed loop
system is in $\mathcal{GPE}$.
\vskip 0.2cm

{\bf C}. If $\Psi\in{\mathcal Odd}$, then, for all $K$ so that
$K+K^*=0$ also the resulting closed loop
system is in ${\mathcal Odd}$.
\end{proposition}

{\bf Proof}\quad
Item {\bf A} is part of \cite[Proposition 8.1 (iii)]{AL2}.
\vskip 0.2cm

{\bf B}. Substituting \eqref{eq:L} in \eqref{Lcl} yields,
for the open loop
\begin{equation}\label{eq:Lol}
A=\left(\begin{smallmatrix}
\hat{A}&~~~\hat{B}\hat{B}^*\\
0&~-\hat{A}^*\end{smallmatrix}\right)\quad\quad
B=\left(\begin{smallmatrix}~~0\\
-\hat{C}^*\end{smallmatrix}\right)\quad\quad
C=(\hat{C}\quad 0)\quad\quad D=0.
\end{equation}
Thus, for the closed loop
\begin{equation}\label{eq:Acl}
A_{\rm c.l.}=\left(\begin{smallmatrix}
\hat{A}&~~~\hat{B}\hat{B}^*\\
-\hat{C}^*K\hat{C}&~-\hat{A}^*\end{smallmatrix}\right).
\end{equation}
Hence,
for all $-K\in\overline{\mathbb P}_p$, this $L_{\rm c.l.}$
satisfies the conditions in \eqref{(ii)d}, so this part of
the claim is established.
\vskip 0.2cm

{\bf C}. Substituting \eqref{RealOdd} in \eqref{Lcl} yields,
\[
A_{\rm c.l.}=\left(\begin{smallmatrix}
T_1+B_1KB_1^*&~~\tilde{A}-B_1KB_2^*\\
\tilde{A}^*+B_2KB_1^*&~~T_2-B_2KB_2^*
\end{smallmatrix}\right)\quad\quad B=
\left(\begin{smallmatrix}B_1\\
B_2\end{smallmatrix}\right)\quad\quad
C=(B_1^*~~~-B_2^*)\quad\quad D=0.
\]
Thus, for all $K$ skew-Hermitian ($K=-K^*$) this $L_{\rm c.l.}$
satisfies condition (ii) in Theorem \ref{main} {\bf B} with
$H={\rm diag}\{-I_{\nu}, I_{n-\nu},~I_p\}$, so the
claim is established
\qed
\vskip 0.2cm

\begin{example}\label{Ex:PseudoSpectFact}
{\rm
We here illustrate item {\bf B} in Proposition
\ref{GPstaticOutputFeedback}. Consider again the
$\mathcal{GPE}$ function $\Psi_2(s)=\frac{-1}{s^2}$
from Example \ref{Ex:GPEImagSpect}. Applying to $\Psi_2(s)$
a static output feedback yields
$\Psi_{2{\rm c.l.}}(s)=\frac{-1}{s^2+k}$ so that
\mbox{$\Psi_{2{\rm c.l.}}\in\mathcal{GPE}$} for all $0>k$.
Indeed, \mbox{$\Psi_{2{\rm c.l.}}(i\omega)=
\frac{1}{\omega^2-k}$} for all $\omega\in\R$.
}
\qed
\end{example}
\vskip 0.2cm

We next use Proposition \ref{GPstaticOutputFeedback} {\bf B}
to specialize the result of Proposition
\ref{Pr:OutputFeedback} to the class of
$\mathcal{GP}$ systems.

\begin{corollary}\label{GPoutputFeedback}
Let $\Psi\in\mathcal{GP}$ be a $p\times p$-valued rational
function s.t. $\lim\limits_{s\rightarrow\infty}\Psi(s)=0$.
Its realization matrix
$L=${\mbox{\tiny$\begin{pmatrix}A&B\\C&0\end{pmatrix}$}}
satisfies \eqref{Lyap} and the closed loop is as in
\eqref{Lcl}. This realization is not minimal, if and only
if there exists \mbox{$\lambda\in{\rm spect}(A)$} so that
\mbox{$\lambda\in{\rm spect}(A_{\rm cl})$}, see \eqref{Lcl},
for all $K$ $-(K+K^*)\in\overline{\mathbb P}_p~$.
\end{corollary}

Indeed, as already pointed out in the proof of part
{\bf B} of Observation \ref{spectL-A}, from
\eqref{eq:symmetry} it follows that $\beta={\rm rank}(B)$
is equal to $\gamma={\rm rank}(C)$. Thus, both directions of
Proposition \ref{Pr:OutputFeedback} hold.
\vskip 0.2cm

Recall that in Theorem \ref{Pr:GPEchar} and the proceeding
discussion we pointed out that \mbox{$\Psi\in\mathcal{GPE}$}
always admits {\em pseudo} spectral factorization, i.e.
\mbox{$\Psi=GG^{\#}$} with $G(s)$ analytic in $\C_-$. It is
desired for $\Psi$ to admit spectral factorization, i.e.
where the factor $G(s)$ is analytic in $\overline{\C_-}$.
Recall also that spectral factorization has various 
applications in control and filtering 
theories, see e.g. \cite{AV}, \cite[Section 16.3]{BGKR}, \cite[Section 6]{FCG},
\cite{Fu}, \cite{HSK}, \cite{KSH}, \cite{LR}, \cite{Ran}.
We now address the following question: Assuming
$\Psi\in\mathcal{GPE}$, does not admit spectral
factorization, under what conditions can one apply
a static output feedback so that
\mbox{$\Psi_{\rm cl}=(I_p-\Psi{K})^{-1}\Psi$}
does admit spectral factorization.
\vskip 0.2cm

Specifically, in Proposition \ref{GPstaticOutputFeedback}
and Corollary \ref{GPoutputFeedback} we examined the use
of static output feedback for moving the poles of
$\Psi\in\mathcal{GP}$ while retaining $\Psi_{\rm cl}$ in
the $\mathcal{GP}$ 
class (and $\mathcal{GPE}$ in particular). We know that a
$\mathcal{GPE}$ function admits spectral factorization if
it is analytic on the imaginary axis. Thus, the problem at
hand is actually about moving by static output feedback
poles of a $\mathcal{GPE}$ function away from the imaginary
axis.
\vskip 0.2cm

To gain intuition we first look at the scalar case. Let
$\Psi\in\mathcal{GPE}$ be written as $\Psi=\frac{N}{D}$
with $N(s)$, $D(s)$ polynomials so that,
$\frac{N(i\omega)}{D(i\omega)}\geq 0$ (including $+\infty$)
for all $\omega\in\R$. Assuming that $N(s)$ and $D(s)$ have
no common roots on $i\R$, this in fact is equivalent to
having for all $\omega\in\R$: $N(i\omega)\geq 0$,
$D(i\omega)\geq 0$ and $D(i\omega)+N(i\omega)>0$. This in
turn is equivalent to having, for all $\omega\in\R$:
$N(i\omega)\geq 0$, $D(i\omega)\geq 0$ and $D(i\omega)+
\alpha{N}(i\omega)>0$ for all $\alpha>0$. Next recall that
\mbox{$\Psi_{\rm cl}=(I_p-\Psi{k})^{-1}\Psi=\frac{N}
{D-kN}$}. In conform with part {\bf B} of
Proposition \ref{GPstaticOutputFeedback},
for all $-k>0$: For all $\omega\in\R$ the closed loop
denominator satisfies
$\left(D(i\omega)-kN(i\omega)\right)>0$. Thus in fact
$\alpha=-k$ and the closed loop denominator is
positive on $i\R$.
To conclude, if $N(s)$ and $D(s)$
have no common imaginary roots\begin{footnote}{In the
spirit of Proposition \ref{Pr:OutputFeedback}, when
restricted to the imaginary axis, the realization is
minimal.}\end{footnote},
for all $-k>0$, $\Psi_{\rm cl}$ is analytic
on $i\R$.
\vskip 0.2cm

Before stating the result we need some preliminaries.
For a pair \mbox{$A\in\C^{n\times n}$},
\mbox{$B\in\C^{n\times m}$} and a region
$\Omega\subseteq\C$, one can define the condition,
\[
{\rm rank}(A-\lambda{I}_n~\vdots~B)=n
\quad\quad\quad\forall~\lambda\in\Omega.
\]
If $\Omega=\C$ the pair $A, B$ is said to be controllable.
If $\Omega$ is a subset of $\C$, this is equivalent to the
existence of a matrix $R\in\C^{m\times n}$ so that
\mbox{${\rm spect}(A+BR)\bigcap\Omega=\emptyset$} (in
particular, for $\Omega=\overline{\C_+}$ the pair is said to
be stabilizable, see e.g. \cite[p. 205]{Ka}
\mbox{\cite[sub-section 4.4]{LR}}).
\vskip 0.2cm

Similarly, for a pair $A\in\C^{n\times n}$,
$C\in\C^{p\times n}$ one can define the condition,
\[
{\rm rank}\left(\begin{smallmatrix}A-\lambda{I}_n\\
\hat{C} \end{smallmatrix}\right)=n\quad\quad\quad
\forall~\lambda\in\Omega.
\]
For $\Omega=\C$ the pair $A, C$ is said to be observable,
for $\Omega=\overline{\C_+}$ the pair is detectable, see
e.g. \mbox{\cite[sub-section 4.4]{LR}}).
\vskip 0.2cm

It is of interest to point here out that in Proposition
\ref{Pr:OutputFeedback} we required minimality of
realization, i.e. $\Omega=\C$. From \cite[Theorem]{KDS}
it follows that having the pair $A, B$ stabilizable and
the pair $A, C$ detectable (i.e. $\Omega=\overline{\C_+}$)
is necessary\begin{footnote}{An additional Riccati type
condition makes it also sufficient.}\end{footnote} for
stabilizability of a system by static output feedback.
Below, we can be ``modest" by resorting to $\Omega=i\R$
and reformulate Proposition \ref{FeedBackSpectFact}.
\vskip 0.2cm

\begin{proposition}\label{FeedBackSpectFact1}
Let $\Psi(s)$ be a $p\times p$-valued rational
$\mathcal{GPE}$ function, which does not admit spectral
factorization. Assume that
$\lim\limits_{s\rightarrow\infty}\Psi(s)=0$ and let
the state space realization be as in \eqref{eq:Lol}.
\vskip 0.2cm

There exists a static output feedback gain $K$,
$-K\in\overline{\mathbb P}_p$ so that $\Psi_{\rm cl}(s)$
admits spectral factorization, if and only if
for $r\in\R$ the two following matrices
\[
\left(\begin{smallmatrix}
\hat{A}-irI_n&\vdots&\hat{B}\end{smallmatrix}\right)
\quad\quad{and}\quad\quad
\left(\begin{smallmatrix}
\hat{A}-irI_n\\ \hat{C}\end{smallmatrix}\right)
\]
%\begin{center}
%{\mbox{\tiny$\begin{pmatrix}
%\hat{A}-irI_n&\vdots&\hat{B}
%\end{pmatrix}$}}
%\quad\quad{and}\quad\quad
%{\mbox{\tiny$\begin{pmatrix}
%\hat{A}-irI_n\\ \hat{C}\end{pmatrix}$}}
%\end{center}
are of full rank.
\end{proposition}
\vskip 0.2cm

{\bf Proof}\quad
First, recall that $\Psi_{\rm cl}(s)$ admits spectral
factorization if and only if in the corresponding
state space realization, the spectrum of
$A_{\rm c.l.}$ avoids the imaginary axis. Following
Proposition \ref{GPstaticOutputFeedback} this in
turn is equivalent to finding conditions on $\hat{A},
\hat{B}$ and $\hat{C}$ so that there exists
\mbox{$-K\in\overline{\mathbb P}_p$} so that
$A_{\rm c.l.}$ in \eqref{eq:Acl}, will have no
eigenvalues on the imaginary axis, i.e. the matrix
\begin{equation}\label{eq:Acl.iR}
A_{\rm c.l.}-irI_{2n}=\left(\begin{smallmatrix}
\hat{A}-irI_n&~~~\hat{B}\hat{B}^*\\
-\hat{C}^*K\hat{C}&~-(\hat{A}-irI_n)^*
\end{smallmatrix}\right)\quad\quad\quad\quad\forall r\in\R
\end{equation}
is nonsingular. Note now that the nonsingularity of
$A_{\rm c.l.}-irI_{2n}$ in \eqref{eq:Acl.iR} is equivalent
to that of
\[
M:=\left(\begin{smallmatrix}0&-I_n\\I_n&~~0\end{smallmatrix}
\right)(A_{\rm c.l.}-irI_{2n})=\left(\begin{smallmatrix}
\hat{C}^*K\hat{C}&~(\hat{A}-irI_n)^*\\
\hat{A}-irI_n&~~~\hat{B}\hat{B}^*\end{smallmatrix}\right).
\]
Namely, we search for $-K\in\overline{\mathbb P}_p$ so that
for all $r\in\R$ the Hermitian matrix $M$ is nonsingular.
\vskip 0.2cm

Next consider the Lyapunov equation $Q=HM+M^*H$ with
\mbox{$H={\rm diag}\{-I_n~,~I_n\}$} and
\[
Q=2{\rm diag}\{-\hat{C}^*K\hat{C},~\hat{B}\hat{B}^*\}.
\]
Now for all $-K\in\overline{\mathbb P}_p$, indeed
$Q\in\overline{\mathbb P}_{2n}$. This implies that the matrix
$M$ has at most $n$ eigenvalues in each open half plane, see
e.g. \cite[Lemma 2.4.5]{HJ2}. Furthermore, from the
Generalized Inertia Theorem for the Lyapunov equation, see
e.g. \cite[Theorems 2.4.7, 2.4.10]{HJ2}, it follows that
there are exactly $n$ eigenvalues in each open half plane,
if and only if there exists $K$, \mbox{$-K\in\overline{\mathbb
P}_p$} so that the pair $M, Q$ is observable.
\vskip 0.2cm

Recall now that the pair $M, Q$ is not
observable, if and only if, there exists a (right)
eigenvector of $M$ which lies in the null-space of $Q$.
Namely, there exists $0\not=v\in\C^{2n}$ so that
\mbox{$Mv=\lambda{v}$} for some $\lambda\in\R$ and
$Qv=0$. This in turn is equivalent to having,
$(M+\frac{1}{2}Q)v=\lambda{v}$ and $Qv=0$. Note
that
\begin{center}
$\hat{M}:=M+\frac{1}{2}Q=${\mbox{\tiny$\begin{pmatrix}
0&~(\hat{A}-irI_n)^*\\
\hat{A}-irI_n&~~~0\end{pmatrix}$}}.
\end{center}
To summarize: The pair $\hat{M}, Q$ is observable if and
only if the pair $M, Q$ is observable and this is
equivalent to the nonsingularity of $M$ (which in turn
is equivalent to the nonsingularity of
$A_{\rm c.l.}-irI_{2n}$).
\vskip 0.2cm

Now, the pair $\hat{M}, Q$ is observable if and only
if there exists $K$, {$-K\in\overline{\mathbb P}_p$}
so that for all $r\in\R$, both matrices:
$\left(\begin{smallmatrix}\hat{A}-irI_n\\
-\hat{C}^*K\hat{C}\end{smallmatrix}\right)$ and
$\left(\begin{smallmatrix}(\hat{A}-irI_n)^*\\
\hat{B}\hat{B}^*\end{smallmatrix}\right)$ are
of full rank.
\vskip 0.2cm

Next, the matrix $\left(\begin{smallmatrix}(\hat{A}-irI_n)^*\\
\hat{B}\hat{B}^*\end{smallmatrix}\right)$ is of full
rank, if and only if the matrix $\left(\hat{A}-irI_n~
\vdots~\hat{B}\right)$ is of full rank, so the first part
of the condition in the claim is established.
\vskip 0.2cm

Note now that if the matrix $\left(\begin{smallmatrix}
\hat{A}-irI_n\\ -\hat{C}^*K\hat{C}\end{smallmatrix}\right)$
is of full rank for some $K$,
{$-K\in\overline{\mathbb P}_p$} it is of full rank
for $K=-I_p$, i.e. the matrix $\left(\begin{smallmatrix}
\hat{A}-irI_n\\ \hat{C}^*\hat{C}\end{smallmatrix}\right)$
is of full rank.  This in turn is equivalent to having
the matrix $\left(\begin{smallmatrix} \hat{A}-irI_n\\
\hat{C}\end{smallmatrix}\right)$ of full rank,
so the second part of the condition
in the claim is established and the proof is complete.
\qed
\vskip 0.2cm

We now illustrate the result of
Proposition \ref{FeedBackSpectFact1}

\begin{example}\label{Ex:SpectFact}
{\rm
Consider the $\mathcal{GPE}$ function $\Psi_2(s)=\frac{-1}{s^2}$
from Examples \ref{Ex:GPEImagSpect} and \ref{Ex:PseudoSpectFact}.
As shown in Example \ref{Ex:GPEImagSpect}. $\Psi_2$ admits only
pseudo-spectral factorization. However, $\Psi_{2{\rm c.l.}}(s)=
\frac{-1}{s^2+k}$ with $0>k$, see Example \ref{Ex:PseudoSpectFact},
can be factored to $\Psi_{2{\rm c.l.}}=G_2G_2^{\#}$ with
\mbox{$G_2(s)=\frac{1}{\sqrt{-k}-s}$}.
From Example \ref{Ex:GPEImagSpect} it follows that
\[
A_{\rm cl}=\left(\begin{smallmatrix}
\hat{A}&~~~\hat{B}\hat{B}^*\\
-k\hat{C}^*\hat{C}&~-\hat{A}^*\end{smallmatrix}\right)=
\left(\begin{smallmatrix}~~0&1\\-k&0\end{smallmatrix}\right)
\quad\quad\quad 0>k.
\]
Indeed, the realization $\hat{A}$, $\hat{B}$, $\hat{C}$
given in Example \ref{Ex:GPEImagSpect} is minimal and
hence in particular the conditions in Proposition
\ref{FeedBackSpectFact1} are satisfied. Thus, indeed the
matrix $A_{\rm cl}$ has no imaginary eigenvalues.
\qed
}
\end{example}
\vskip 0.2cm

From the proof of Proposition \ref{FeedBackSpectFact1}
it follows that if $\Psi\in\mathcal{GPE}$ and $K$,
$-K\in\overline{\mathbb P}_p$
are so that \mbox{$\Psi_{\rm cl}=(I_p-\Psi{K})^{-1}\Psi$}
admits spectral factorization, the same is true for
$\alpha{K}$ where $\alpha>0$ may be arbitrarily small. Namely,
$\Psi_{\rm cl}(s)$ may be a small perturbation of $\Psi(s)$.
\vskip 0.2cm

We conclude this section by noting that upon comparing
Propositions \ref{Pr:OutputFeedback} and
\ref{FeedBackSpectFact1} one can make the following
statement.

\begin{corollary}
Let $G(s)=\hat{C}(sI_n-\hat{A})^{-1})\hat{B}$ be a 
realization of a rational function so that both
$\hat{B}, \hat{C}\in\C^{p\times n}$ are of the same
rank. The following are equivalent.
\begin{itemize}
\item{}There exists a static output feedback gain
$K\in\C^{p\times p}$ so that the closed loop system
\mbox{$(I_p-GK)^{-1}G$} is analytic on $i\R$.

\item{}There exists a static output feedback with
$-K\in\overline{\mathbb P}_p$  so that the closed loop
system $(I_p-GG^{\#}K)^{-1}GG^{\#}$ is analytic on
$i\R$.
\end{itemize}
\end{corollary}

\end{document}